\documentclass[a4paper,12pt]{amsart}
\usepackage{amssymb}
\usepackage{graphicx}
\usepackage{ifthen}
\nonstopmode \numberwithin{equation}{section}
\newtheorem{thm}{Theorem}%[section]
\newtheorem{cor}{Corollary}%[section]
\newtheorem{lem}{Lemma}%[section]
\newtheorem{conj}{Conjecture}

\theoremstyle{definition}
\newtheorem{defn}{Definition}[section]

\newtheorem{prob}[equation]{Problem}

\newenvironment{rem}{%
\bigskip
\noindent \textsl{{\sl Remark. }}}{\bigskip}
\newenvironment{rems}{%
\bigskip
\noindent \textsl{{\sl Remarks. }}}{\bigskip}

\newcounter {own}
\def\theown {\thesection       .\arabic{own}}

\newenvironment{pf}[1][]{%
 \vskip 3mm
 \noindent
 \ifthenelse{\equal{#1}{}}%
  {{\slshape Proof. }}%
  {{\slshape #1.} }%
 }%
{\qed\bigskip}

\newcounter{alphabet}
\newcounter{tmp}
\newenvironment{Thm}[1][]{\refstepcounter{alphabet}%
\bigskip%
\noindent%
{\bf Theorem \Alph{alphabet}}%
\ifthenelse{\equal{#1}{}}{}{ (#1)}%
{\bf .} \itshape}{\vskip 8pt}

\newcommand{\A}{{\mathcal A}}
\newcommand{\B}{{\mathcal B}}

\newcommand{\U}{{\mathcal U}}
\newcommand{\es}{{\mathcal S}}

\newcommand{\ID}{{\mathbb D}}
\newcommand{\IN}{{\mathbb N}}
\newcommand{\IC}{{\mathbb C}}

\newcommand{\CC}{{\mathcal C}}

\def\be{\begin{equation}}
\def\ee{\end{equation}}

\newcommand{\bee}{\begin{enumerate}}
\newcommand{\eee}{\end{enumerate}}

\newcommand{\blem}{\begin{lem}}
\newcommand{\elem}{\end{lem}}
\newcommand{\bthm}{\begin{thm}}
\newcommand{\ethm}{\end{thm}}
\newcommand{\bcor}{\begin{cor}}
\newcommand{\ecor}{\end{cor}}
\newcommand{\beg}{\begin{examp}}
\newcommand{\eeg}{\end{examp}}
\newcommand{\begs}{\begin{examples}}
\newcommand{\eegs}{\end{examples}}
\newcommand{\bdefe}{\begin{defn}}
\newcommand{\edefe}{\end{defn}}
\newcommand{\bprob}{\begin{prob}}
\newcommand{\eprob}{\end{prob}}
\newcommand{\bei}{\begin{itemize}}
\newcommand{\eei}{\end{itemize}}

\newcommand{\bcon}{\begin{conj}}
\newcommand{\econ}{\end{conj}}
\newcommand{\bcons}{\begin{conjs}}
\newcommand{\econs}{\end{conjs}}
\newcommand{\bprop}{\begin{propo}}
\newcommand{\eprop}{\end{propo}}
\newcommand{\br}{\begin{rem}}
\newcommand{\er}{\end{rem}}
\newcommand{\brs}{\begin{rems}}
\newcommand{\ers}{\end{rems}}
\newcommand{\bo}{\begin{obser}}
\newcommand{\eo}{\end{obser}}
\newcommand{\bos}{\begin{obsers}}
\newcommand{\eos}{\end{obsers}}
\newcommand{\bpf}{\begin{pf}}
\newcommand{\epf}{\end{pf}}
\newcommand{\ba}{\begin{array}}
\newcommand{\ea}{\end{array}}
\newcommand{\beq}{\begin{eqnarray}}
\newcommand{\beqq}{\begin{eqnarray*}}
\newcommand{\eeq}{\end{eqnarray}}
\newcommand{\eeqq}{\end{eqnarray*}}

\newcommand{\ov}{\overline}

\newcounter{minutes}
\divide\time by 60
\newcounter{hours}
\multiply\time by 60 \addtocounter{minutes}{-\time}
\begin{document}
\title{Bohr phenomenon for locally univalent functions and logarithmic power series}
\begin{center}
{\tiny \texttt{FILE:~\jobname .tex,
        printed: \number\year-\number\month-\number\day,
        \thehours.\ifnum\theminutes<10{0}\fi\theminutes}
}
\end{center}
\author{Bappaditya Bhowmik${}^{\mathbf{*}}$}
\address{Bappaditya Bhowmik, Department of Mathematics,
Indian Institute of Technology Kharagpur, Kharagpur - 721302, India.}
\email{bappaditya@maths.iitkgp.ac.in}
\author{Nilanjan Das}
\address{Nilanjan Das, Department of Mathematics,
Indian Institute of Technology Kharagpur, Kharagpur - 721302, India.}
\email{nilanjan@iitkgp.ac.in}

\subjclass[2010]{30B10, 30C60, 31A05, 30C45}
\keywords{Bohr radius, Locally univalent functions, Logarithmic coefficients.\newline
${}^{\mathbf{*}}$ Corresponding author}
%\date{ %\today Feb. 5, 08
%May 18, 2017; File: Bound-BB.tex}
%\date{%\today
%September 01, 06; File: lau$_{-}$revise1.tex}

\begin{abstract}
In this article we prove Bohr inequalities for sense-preserving $K$-quasiconformal harmonic mappings defined in $\ID$
and obtain the corresponding results for sense-preserving harmonic mappings by letting $K\to\infty$. One of the results
includes the sharpened version of a theorem by Kayumov \emph{et. al.} (\textit{Math. Nachr.}, 291 (2018), no. 11--12, 1757--1768).
In addition Bohr inequalities have been established for uniformly locally univalent holomorphic functions,
and for $\log(f(z)/z)$ where $f$ is univalent or inverse of a univalent function.
\end{abstract}
\thanks{The first author of this article would like to thank
SERB, DST, India (Ref.No.- MTR/2018/001176) for its financial support through MATRICS grant.}

\maketitle
\pagestyle{myheadings}
\markboth{B. Bhowmik, N. Das}{Bohr phenomenon for locally univalent functions and logarithmic power series}

\bigskip
%\noindent
\section{Introduction}
The origin of Bohr phenomenon lies in the seminal work by Harald Bohr \cite{Bohr}, which included the following
(improved) result.

\begin{Thm}\label{P1TheoA}
Let $f(z)= \sum_{n=0}^{\infty}a_{n}z^n$ be holomorphic in the open unit disk $\ID$ and $|f(z)|<1$ for all $z\in\ID$, then
\be\label{P3eq27}
\sum_{n=0}^{\infty}|a_n|r^n\leq 1
\ee
for all $z\in\ID$ with $|z|=r\leq1/3$.
\end{Thm}
%The constant $1/6$ was subsequently sharpened to $1/3$ by Wiener, Riesz and Schur.
Inequalities of similar nature are
being extensively investigated nowadays in different frameworks, and have become famous by the name \emph{Bohr inequalities}.
To have a glimpse of the ongoing current research in Bohr radius problem the reader is urged to glance through some of the recent articles,
e.g. \cite{Aiz1, Boas, Def1, Paul3} and the references therein.
Now we concentrate on a generalized treatment of the Bohr radius problem introduced in \cite{Abu}, using the concept of subordination.
For two holomorphic functions $f$ and $g$ in $\ID$,
we say $g$ is subordinate to $f$ if there exists a function $\phi$, holomorphic
in $\ID$ with $\phi(0)=0$ and $|\phi(z)|<1$, satisfying $g=f\circ\phi$.
Throughout this article we denote $g$ is subordinate to $f$  by $g\prec f$. Also
the class of functions $g$ subordinate to a fixed function $f$ will be denoted by $S(f)$. Now according to \cite{Abu}
we say that $S(f)$ has Bohr phenomenon if for any $g(z)=\sum_{n=0}^{\infty}b_{n}z^n\in S(f)$, there is a $r_0\in(0,1]$
such that
\be\label{P3eq28}
\sum_{n=1}^{\infty}|b_{n}|r^n\leq d(f(0),\partial f(\ID))
\ee
for $|z|=r<r_0$. Here $d(f(0),\partial f(\ID))$ denotes the Euclidean distance between $f(0)$ and the boundary of domain $f(\ID)$.
It is seen that whenever a holomorphic function $g$ maps $\ID$ into a domain $\Omega$ other than $\ID$, then in a general sense the
Bohr inequality $(\ref{P3eq28})$ can be established if $g$ can be recognized as a member of $S(f)$, $f$ being the covering map
from $\ID$ onto $\Omega$ satisfying $f(0)=g(0)$.
In particular, if we take $\Omega=\ID$, then for any holomorphic $g:\ID\to\Omega$ there exists a disk automorphism $f$
such that $g(0)=f(0)$ and $g\in S(f)$.
In this case $ d(f(0),\partial\ID)=1-|f(0)|$, and hence $(\ref{P3eq28})$ reduces to $(\ref{P3eq27})$. Bohr phenomenon has been explored
using the above definition in a number of papers, e.g. \cite{Abu, Abu1, Abu3, Bar, BB}.
One of the goals of the present article is to extend the Bohr inequalities of type $(\ref{P3eq28})$
for certain harmonic functions in a suitable fashion.
%We now present some definitions which are required for our further discussions.
A complex valued function $f(z)=u(x,y)+iv(x,y)$ of
$z=x+iy\in\ID$ is called \emph{harmonic} if both $u$ and $v$ satisfy the Laplace's equation
$$
\frac{\partial^2H}{\partial x^2}+ \frac{\partial^2H}{\partial y^2}=0,
$$
where $H(x,y)$ is a real valued function. It is well known that under the assumption $g(0)=0$,
$f$ has a unique \emph{canonical representation} $f=h+\ov{g}$,
where $h$ and $g$ are holomorphic in $\ID$. In view of this representation, $f$ is locally univalent and sense-preserving
whenever the Jacobian $J_f(z):=|h^{\prime}(z)|^2-|g^{\prime}(z)|^2>0$ for all $z\in\ID$. A sense-preserving homeomorphism
defined in $\ID$ which is also harmonic is called $K$-\emph{quasiconformal}, $K\in[1,\infty)$ if
the (second complex) dilatation $w_f:=g^\prime/h^\prime$ satisfies
$|w_f(z)|\leq k$, $k=(K-1)/(K+1)\in[0,1)$. Now it is easy to see that the
aforesaid definitions and notations for subordination of holomorphic functions
can be adopted for harmonic functions without any change (cf. \cite{Lis}).
In the present day theory of harmonic mappings, investigations are often carried out to explore the connections
between the holomorphic part, or some suitable holomorphic counterpart of a given harmonic mapping and the map itself
(see f.i. \cite{Her, Par}). Motivated by this perspective, in this article we prove Bohr inequalities similar to
$(\ref{P3eq28})$ for $S(f)$ under the assumption that $f$ is
a sense-preserving $K$-quasiconformal harmonic mapping  defined in $\ID$, where the holomorphic part
$h$ is univalent or convex univalent.
%In Theorem 1 and Remark 1 from \cite{Abu}, Bohr radii for $S(f)$ were obtained for univalent
%and convex univalent holomorphic function $f$ respectively.
Further, as another application of the technique used in proving this
theorem, we establish the sharpened version of \cite[Theorem 3.1]{Pon1}.
We here mention that a number of Bohr inequalities
for sense-preserving $K$-quasiconformal harmonic mappings
have been obtained in \cite{Pon1}, which mostly bear the classical flavor of the Bohr radius problem.

We now turn our attention to the class $\mathcal{H}$ of complex valued holomorphic functions $f$ defined in $\ID$.
A considerably interesting subfamily of $\mathcal{H}$ is the class of
uniformly locally univalent functions (see \cite{Kim, Pom, Su1, Ya}). Here we clarify that
$f\in\mathcal{H}$ is said to be \textit{uniformly locally univalent} if
there exists $a>0$ such that $f$ is univalent
on each hyperbolic disk $\ID_a^h(z_0):=\{z\in\ID: |(z-z_0)/(1-\ov{z_0}z)|<\tanh a\}$ with center $z_0\in\ID$
and radius $a$.
It is well known that (cf. \cite{Kim, Ya}) a function $f\in\mathcal{H}$ is uniformly locally univalent if and only if the
pre-Schwarzian norm
$$
\|P_f\|:=\sup_{z\in\ID}(1-|z|^2)|f^{\prime\prime}(z)/f^\prime(z)|<\infty.
$$
Now let $\A:=\{f\in\mathcal{H}:f(0)=f^\prime(0)-1=0\}$. Since $f^{\prime\prime}/f^\prime$ remains invariant under
the post-composition by a non-constant linear function,
in view of the above characterization it is quite natural to consider the class
$\B(\lambda):=\{f\in\A: \|P_f\|\leq 2\lambda\}$ for any $\lambda\in[0,\infty)$ (compare \cite{Kim}).
In this paper we derive a Bohr inequality of type $(\ref{P3eq28})$ for the functions in $\B(\lambda)$. As
$\B(0)=\{z\}$, we consider $\lambda\in(0,\infty)$ only to prove the result.

Before we proceed further we need to introduce the following subfamilies of $\A$ to facilitate our discussion.
Let the subclass of univalent functions in $\A$ be denoted by $\es$.
Two well known subclasses of $\es$ are $\es^*$ and $\CC$ which consist of starlike and convex univalent functions respectively.
We now consider the \emph{logarithmic coefficients} of any $f\in\es$, which are defined by
\be\label{P3eq2}
\log \left(f(z)/z\right)=2\sum_{n=1}^\infty \gamma_nz^n, z\in\ID
\ee
(see, f.i. \cite[p. 151]{Dur}). The importance of logarithmic coefficients in univalent function theory is already well regarded
due to the substantial role played by them in the proof of Bieberbach conjecture.
We know that proving inequalities concerning $|\gamma_n|$'s is considered
to be a challenging problem till date (see f.i. \cite{Pon3, Roth} and references therein) due to the unavailability
of the sharp bounds on $|\gamma_n|$'s for $n\geq3$, where $f\in\es$.
Inspired by this fact, in this article we have considered the problem of establishing Bohr inequalities
similar to the inequality $(\ref{P3eq27})$ for $\log(f(z)/z)$. More
precisely, we will say that $\log(f(z)/z)$ has Bohr radius $r_0\in(0,1]$ if
\be\label{P3eq3}
2\sum_{n=1}^\infty|\gamma_n|r^n\leq 1
\ee
for $|z|=r< r_0$.
We comment here that the quantity $d(\log(f(z)/z), \partial\Omega)$,
where $\Omega$ is the image of $\ID$ under the function $\log(f(z)/z)$, can be an arbitrarily small
positive number for $f\in\es$. One can easily see this by choosing the univalent polynomials $f_n(z)=z+(z^2/n), z\in\ID$ for each $n\geq 2$
(cf. \cite[p. 267]{Dur}), and observing that the image of $\log(f_n(z)/z)$ does not include the point $\log(1+(1/n))$.
This fact backs up our choice to define the Bohr phenomenon
for $\log(f(z)/z)$ in the classical manner instead of using any inequality of the type $(\ref{P3eq28})$.
We derive Bohr inequalities in the form of $(\ref{P3eq3})$ while $f$ is a member of $\es, \es^*, \CC$.
Moreover, for $f\in\es(\mbox{or}\, \es^*)$, $f^{-1}(w)$ is defined in a neighborhood of the origin, which in particular can be chosen to be
$\ID_{1/4}:=\{w\in\IC:|w|<1/4\}$, as we know that any $f\in\es(\mbox{or}\, \es^*)$ covers $\ID_{1/4}$. Therefore it is possible to define the
logarithmic coefficients of $f^{-1}$ for $f\in\es(\mbox{or}\,\es^*)$ by the following expression:
\be\label{P3eq11}
\log \left(f^{-1}(w)/w\right)=2\sum_{n=1}^\infty \gamma_nw^n,\, w\in\ID_{1/4}
\ee
(compare \cite{Pon}). We also compute the Bohr radius for $\log \left(f^{-1}(w)/w\right)$ with respect to the
inequality $(\ref{P3eq3})$, where $r=|w|$ and $f\in\es(\mbox{or}\,\es^*)$.
Another important class $\U(\lambda)$ is being extensively studied by many authors (cf. \cite{Obra4, Pon3} and the references therein)
which is defined by $\U(\lambda)=\{f\in\A: |U_f(z)|<\lambda\}$
where $0<\lambda\leq 1$, and
$$
U_f(z):=\left(z/f(z)\right)^2f^\prime(z)-1, z\in\ID.
$$
It is well known that $\U(\lambda)\subsetneq\es$, and also that $\U(\lambda)$ neither contains
$\es^*$ nor is contained in it. Since the coefficient problem for $f$ or $\log(f(z)/z), f\in\U(\lambda)$
has not yet been fully solved, the Bohr radius problem for $\log(f(z)/z)$ becomes
quite appealing whenever $f\in\U(\lambda)$. Therefore, we end this article with a Bohr inequality for $\log(f(z)/z)$, $f\in\U(\lambda)$.
Here we mention that the power series described in $(\ref{P3eq2})$ and $(\ref{P3eq11})$ will be called \textit{logarithmic power series}
in this article.
\section{Bohr phenomenon for locally univalent functions}
We prove the following lemma which will be required to establish  next two  theorems in this section.
\blem\label{P3lem1}
Let $h(z)=\sum_{n=0}^\infty a_nz^n$ and $g(z)=\sum_{n=0}^\infty b_nz^n$ be two holomorphic functions defined in $\ID$ such that
$g(z)=M\phi(z)h(z)$ for some $M>0$, $\phi:\ID\to\ID$ being a holomorphic function with an expansion
$\phi(z)=\sum_{n=0}^\infty c_nz^n$. Then
\be\label{P3eq18}
\sum_{n=0}^\infty |b_n|r^n\leq M\sum_{n=0}^\infty |a_n|r^n
\ee
for $|z|=r\leq 1/3$.
\elem
\bpf
Since $g(z)=M\phi(z)h(z)$, taking the Cauchy product of two series, we get
$$
\sum_{n=0}^\infty b_nz^n=M\left(\sum_{i=0}^\infty c_iz^i\right)\left(\sum_{j=0}^\infty a_jz^j\right)
=M\sum_{n=0}^\infty\left(\sum_{t=0}^na_tc_{n-t}\right)z^n.
$$
Hence for $|z|=r\in[0,1)$,
\be\label{P3eq16}
\sum_{n=0}^\infty |b_n|r^n=M\sum_{n=0}^\infty\left|\sum_{t=0}^na_tc_{n-t}\right|r^n
\leq M\sum_{n=0}^\infty\left(\sum_{t=0}^n|a_t||c_{n-t}|\right)r^n.
\ee
Now it is easy to see that
\be\label{P3eq17}
M\sum_{n=0}^\infty\left(\sum_{t=0}^n|a_t||c_{n-t}|\right)r^n
=M\left(\sum_{i=0}^\infty |c_i|r^i\right)\left(\sum_{j=0}^\infty |a_j|r^j\right).
\ee
From the Theorem~A it is known that $\sum_{i=0}^\infty |c_i|r^i\leq 1$ for $r\leq 1/3$.
Therefore combining $(\ref{P3eq16})$ and $(\ref{P3eq17})$, we get that the inequality $(\ref{P3eq18})$ holds
for $|z|=r\leq 1/3$.
\epf

We are now ready to prove the first theorem of this section,
which includes sharp Bohr radius for the subordinating family of a sense-preserving
$K$-quasiconformal harmonic mapping with univalent holomorphic part.
\bthm\label{P3thm7}
Let $f(z)=h(z)+\ov{g(z)}=\sum_{n=0}^\infty a_nz^n+\ov{\sum_{n=1}^\infty b_nz^n}$ be a sense-preserving $K$-quasiconformal harmonic mapping
defined in $\ID$ such that $h$ is univalent in $\ID$, and
let $f_1(z)=h_1(z)+\ov{g_1(z)}=\sum_{n=0}^\infty c_nz^n+\ov{\sum_{n=1}^\infty d_nz^n}\in S(f)$.
Then
\be\label{P3eq24}
\sum_{n=1}^\infty |c_n|r^n+\sum_{n=1}^\infty |d_n|r^n\leq d(h(0), \partial h(\ID))
\ee
for $|z|=r\leq r_0=(5K+1-\sqrt{8K(3K+1)})/(K+1)$. This result is sharp for the function $p(z)=z/(1-z)^2+k\ov{z/(1-z)^2}$,
where $k=(K-1)/(K+1)$. Moreover if we take $h$ to be convex univalent then the inequality $(\ref{P3eq24})$ holds for
$r\leq r_0=(K+1)/(5K+1)$. This result is again sharp for the function $q(z)=z/(1-z)+k\ov{z/(1-z)}$.
\ethm
\bpf
From the definition of sense-preserving $K$-quasiconformal harmonic mappings, $h^\prime(z)\neq 0$ for all $z\in\ID$, and
the dilatation $w_f=g^\prime/h^\prime$ satisfies $|w_f(z)|\leq k<1, z\in\ID$ where $k=(K-1)/(K+1)$.
From maximum modulus principle, if $w_f$ is non-constant then $|w_f(z)|< k$ for all $z\in\ID$.
Therefore assuming $w_f$ non-constant, we see that there exists
a holomorphic function $\phi:\ID\to\ID$ such that
$g^\prime(z)=k\phi(z)h^\prime(z), z\in\ID$. An application of Lemma \ref{P3lem1} readily gives
$$
\sum_{n=1}^\infty n|b_n|r^{n-1}\leq k\sum_{n=1}^\infty n|a_n|r^{n-1}
$$
for $r\leq 1/3$, which, upon integration from $0$ to $r$ gives
\be\label{P3eq25}
\sum_{n=1}^\infty|b_n|r^n\leq k\sum_{n=1}^\infty|a_n|r^n
\ee
for $r\leq 1/3$.
Now it is well known that since $h$ is univalent,
$|a_1|\leq4d(h(0), \partial h(\ID))$ (see f.i. \cite[Lemma 1]{Abu}), and the famous de Branges's
theorem asserts that $|a_n|\leq n|a_1|$ for $n\geq 1$. Consequently $|a_n|\leq 4nd(h(0), \partial h(\ID))$ for
all $n\geq 1$. Therefore from $(\ref{P3eq25})$ we get
\be\label{P3eq26}
\sum_{n=1}^\infty |a_n|r^n+\sum_{n=1}^\infty |b_n|r^n\leq (1+k)\sum_{n=1}^\infty|a_n|r^n
\leq \frac{4(1+k)r}{(1-r)^2}d(h(0), \partial h(\ID))
\ee
for $r\leq 1/3$. From a direct computation we obtain that the right hand side
of the inequality $(\ref{P3eq26})$ is less or equal to $d(h(0), \partial h(\ID))$ if
$r^2-(6+4k)r+1\geq 0$, or equivalently if $r\leq r_0=(5K+1-\sqrt{8K(3K+1)})/(K+1)$.
Again by straightforward calculations one can verify that $r_0<1/3$.
Therefore to prove our first assertion in the theorem,
it suffices to show that for $r\leq 1/3$
$$
\sum_{n=1}^\infty |c_n|r^n+\sum_{n=1}^\infty |d_n|r^n\leq\sum_{n=1}^\infty |a_n|r^n+\sum_{n=1}^\infty |b_n|r^n,
$$
which is indeed true, because $f_1\prec f$ implies $h_1\prec h$ and $g_1\prec g$ (cf. \cite[p. 164, Sec. 2]{Lis}), and
therefore from \cite[Lemma 1]{BB}, $\sum_{n=1}^\infty |c_n|r^n\leq \sum_{n=1}^\infty |a_n|r^n$ and
$\sum_{n=1}^\infty |d_n|r^n\leq \sum_{n=1}^\infty |b_n|r^n$ respectively for $r\leq 1/3$.
Now if $w_f$ is constant, then $w_f=ck$ for some $|c|=1$, i.e. $g^\prime(z)=ckh^\prime(z)$.
As a result equality occurs in $(\ref{P3eq25})$ for all $r<1$, and hence this case can be settled by following the
similar lines of reasoning we have already used.
The sharpness part for the function $p$ can be verified from direct calculations.

If $h$ is taken to be convex univalent, we only need to note that $|a_n|\leq |a_1|$ for $n\geq 1$ and
$|a_1|\leq 2d(h(0), \partial h(\ID))$ (see, for example \cite[Lemma 2]{Abu}). Rest of the proof can be completed
by following similar lines of argument presented above.
\epf

\brs
In connection with the above theorem the following interesting observations are made.
\bee
\item
The Theorem 1 and Remark 1 from \cite{Abu} are special instances of the above Theorem \ref{P3thm7}, obtained by setting
$K=1$.
\item
Letting $K\to\infty$ we get that $(\ref{P3eq24})$ holds for $r\leq r_0=5-2\sqrt{6}$, where $f$ is a sense-preserving harmonic mapping
defined in $\ID$ with $h$ univalent, and for $r\leq r_0=1/5$ with $h$ convex univalent. Both of these radii are the best possible.
\eee
\ers

In the next theorem we prove sharp Bohr inequality for a sense-preserving $K$-quasiconformal harmonic mapping $f$
with the canonical representation $f=h+\ov{g}$, under the additional assumptions that $h$ is bounded and $g^\prime(0)=0$.

%and its corollary refine the
%Theorem $3.1$ and Corollary $3.2$ respectively from the recent article \cite{Pon1}.
%The reader is urged to glance through \cite{Lis} for the definition and properties of subordination on harmonic functions.
\bthm\label{P3thm5}
Let $f(z)=h(z)+\ov{g(z)}=\sum_{n=0}^\infty a_nz^n+\ov{\sum_{n=2}^\infty b_nz^n}$ be a sense-preserving $K$-quasiconformal harmonic mapping
defined in $\ID$, where $h$ is bounded on $\ID$. Then
\be\label{P3eq13}
\sum_{n=0}^\infty |a_n|r^n+\sum_{n=2}^\infty |b_n|r^n\leq\|h\|_\infty
\ee
for $|z|=r\leq r_0$, where $r_0$ is the only root in $(0,1)$ of the equation
\be\label{P3eq14}
\frac{4Kr}{(K+1)(1-r)}+\frac{2(K-1)\log(1-r)}{K+1}=1.
\ee
This $r_0$ is the best possible.
\ethm
\bpf
Without loss of generality we can consider $\|h\|_\infty=1$. Also we observe that the case $K=1$
follows from Theorem~A. Hence it is enough to consider $K>1$.
As in the proof of Theorem \ref{P3thm7}, it is easily seen that if $w_f\not\equiv 0$, then
$w_f/k$ is a holomorphic self mapping of $\ID$ with $w_f(0)=0$. From the
Schwarz lemma we can conclude that $\phi(z):=w_f(z)/kz$ is again a holomorphic self mapping of $\ID$
if $w_f(z)\neq kcz$, $z\in\ID$ for some $|c|=1$.
%Now let $\phi$ has a Taylor expansion
%\be\label{P3eq15}
%\phi(z)=\sum_{n=0}^\infty c_nz^n.
%\ee
Therefore assuming $w_f(z)\neq kcz$, a use of Lemma~\ref{P3lem1} on $g^\prime(z)=kz\phi(z)h^\prime(z)$ yields
$$
\sum_{n=2}^\infty n|b_n|r^{n-1}\leq k\sum_{n=1}^\infty n|a_n|r^n,
$$
for $|z|=r\leq 1/3$, or equivalently
$$
\sum_{n=1}^\infty(n+1)|b_{n+1}|r^n\leq k\sum_{n=1}^\infty n|a_n|r^n
$$
for $r\leq 1/3$. Integrating both sides of the above inequality from $0$ to $r$ we have, for $r\leq 1/3$:
$$
\sum_{n=1}^\infty|b_{n+1}|r^{n+1}\leq k\sum_{n=1}^\infty\left(\frac{n}{n+1}\right)|a_n|r^{n+1},
$$
which is same as saying
\be\label{P3eq19}
\sum_{n=2}^\infty|b_{n}|r^{n}\leq k\sum_{n=2}^\infty\left(\frac{n-1}{n}\right)|a_{n-1}|r^{n}.
\ee
Therefore using $(\ref{P3eq19})$ and the well known estimates $|a_n|\leq 1-|a_0|^2$ for $n\geq 1$, we have,
for $r\leq 1/3$:
\be\label{P3eq20}
\sum_{n=0}^\infty |a_n|r^n+\sum_{n=2}^\infty |b_n|r^n
\leq |a_0|+(1-|a_0|^2)\sum_{n=1}^\infty r^n+k(1-|a_0|^2)\sum_{n=2}^\infty\left(\frac{n-1}{n}\right)r^n.
\ee
We here mention that for  the case $w_f(z)=kcz$ , $|c|=1$ the inequality $(\ref{P3eq20})$ can be obtained
from direct calculation, i.e. without any use of Lemma \ref{P3lem1}, and will hold for all $r<1$.
A little computation will now reveal that the right hand side of $(\ref{P3eq20})$ is equal to
$|a_0|+(1-|a_0|^2)\left((1+k)(r/(1-r))+k\log(1-r)\right),$
which is less or equal to $1$ if
$$
(1+|a_0|)\left(\frac{2Kr}{(K+1)(1-r)}+\frac{(K-1)\log(1-r)}{K+1}\right)\leq 1,
$$
which is again true if
$$
\psi(r):=\frac{4Kr}{(K+1)(1-r)}+\frac{2(K-1)\log(1-r)}{K+1}-1\leq0.
$$
To prove the theorem it is sufficient to show that $\psi(r)$ has exactly one root $r_0$ in $(0,1)$,
$r_0< 1/3$ and $\psi(r)\leq 0$ if and only if $r\leq r_0$. We observe that $\psi(0)=-1<0$ and
$\psi(1/3)=(K-1)(1+\log4-\log 9)/(K+1)>0$. By intermediate value property of continuous functions
there exists $r_0\in(0, 1/3)$ such that $\psi(r_0)=0$. Moreover, we observe that for all $r\in(0,1)$,
%$$
%\psi^\prime(r)= \frac{2\left(1+\frac{K(1+r)}{1-r}\right)}{(1-r)(1+K)}>0,
%$$
$\psi^\prime(r)>0$; which implies $\psi$ is strictly increasing in $(0,1)$. This asserts that $r_0$ is the only root of
$\psi$ in $(0,1)$, and that $\psi(r)\leq 0\Leftrightarrow r\leq r_0$. To see that $r_0$ is best possible
one can refer to the computations from \cite[p. 1763]{Pon1} included in the proof of \cite[Theorem 3.1]{Pon1}.
\epf
\bcor\label{P3cor1}
Let $f(z)=h(z)+\ov{g(z)}=\sum_{n=0}^\infty a_nz^n+\ov{\sum_{n=2}^\infty b_nz^n}$ be a sense-preserving harmonic mapping
defined in $\ID$, where $h$ is bounded on $\ID$. Then inequality $(\ref{P3eq13})$ holds for
$|z|=r\leq r_0=0.299\cdots$, where $r_0$ is the only root in $(0,1)$ of the equation
%\be\label{P3eq21}
$$
\frac{4r}{1-r}+2\log(1-r)=1.
$$
%\ee
This $r_0$ is the best possible.
\ecor
\bpf
Follows immediately by letting $K\to\infty$ in the equation $(\ref{P3eq14})$.
\epf

\brs
The following comments are in order.
\bee
\item
The Theorem \ref{P3thm5} (resp. Corollary \ref{P3cor1}) is the refined version of
\cite[Theorem 3.1]{Pon1} (resp. \cite[Corollary 3.2]{Pon1}). It may be observed that
our proof uses a substantially different method compared to \cite[Theorem 3.1]{Pon1}.
Also Corollary \ref{P3cor1} can be taken as a generalization of Proposition 1 from \cite[p. 211]{Pon2}, provided we note that
the conclusion of Corollary \ref{P3cor1} remains unchanged if the hypothesis ``sense-preserving" is replaced by
$|g^\prime(z)|\leq|zh^\prime(z)|$, $z\in\ID$.
\item
It is interesting to note that Theorem $1.1$ and Theorem $1.3$ (and hence Corollary $1.4$ and Corollary $1.5$) from the paper
\cite{Pon1} can also be established using the Lemma~\ref{P3lem1}. One should, however, note that the part of
Corollary 1.4 which remarks on the cases that $a_0= 0$ or $|a_0|$ being replaced by $|a_0|^2$ would produce a better
Bohr radius $1/3$ instead of $1/5$, has to be proved separately, as \cite[Theorem 1.2]{Pon1} can not be derived from the Lemma 1.
\eee
\ers

We now establish a (possibly non-sharp) Bohr inequality for normalized uniformly locally univalent holomorphic
functions with bounded pre-Schwarzian norm.
\bthm\label{P3thm6}
Let $f(z)=z+\sum_{n=2}^\infty a_nz^n\in\B(\lambda)$ for some $\lambda\in(0,\infty)$. Then the inequality
\be\label{P3eq22}
r+\sum_{n=2}^\infty|a_n|r^n\leq d(f(0), \partial f(\ID))
\ee
holds for $|z|=r\leq r_0$, where $r_0$ is the only root in $(0,1)$ of the equation
\be\label{P3eq23}
r+r\sqrt{\exp\left(4\lambda^2r^2/(1-r^2)\right)-1}\sqrt{(\pi^2/6)-1}=-F_\lambda(-1).
\ee
The function $F_\lambda$ is given by
$$
F_\lambda(z):=\int_0^z\left(\frac{1+t}{1-t}\right)^\lambda dt.
$$
\ethm
\bpf
From the definition of $\B(\lambda)$ we see, for any $|z|=r$
$$
|(\log f^\prime(z))^\prime|\leq \frac{2\lambda}{1-r^2},
$$
where $\log f^\prime(z)=\sum_{n=1}^\infty c_nz^n, z\in\ID$. Now the above inequality implies
$$
\sum_{n=1}^\infty n|c_n|^2r^{2n}=\frac{1}{\pi}\iint\limits_{|\xi|<r}
|(\log f^\prime(\xi))^\prime|^2R\,dR\,d\theta\leq
\frac{4\lambda^2}{\pi}\iint\limits_{|\xi|<r}\frac{R}{(1-R^2)^2}dR\,d\theta=\frac{4\lambda^2r^2}{1-r^2}.
$$
We clarify that the dummy variable $\xi$ inside the integration is taken to be $\xi=Re^{i\theta}$,
$0\leq R\leq r$ and $0\leq\theta<2\pi$.
Using a minor variant of first Lebedev-Milin inequality (see \cite[pp. 143-144]{Dur}) we obtain
$$
\sum_{n=1}^\infty n^2|a_n|^2r^{2n-2}\leq \exp\left(\sum_{n=1}^\infty n|c_n|^2r^{2n}\right)\leq
\exp\left(\frac{4\lambda^2r^2}{1-r^2}\right).
$$
In other words
$$
\sum_{n=2}^\infty n^2|a_n|^2r^{2n}\leq r^2\left(\exp\left(\frac{4\lambda^2r^2}{1-r^2}\right)-1\right),
$$
which, along with an application of Cauchy-Schwarz inequality yields
$$
\sum_{n=1}^\infty|a_n|r^n\leq r+ \sqrt{\sum_{n=2}^\infty n^2|a_n|^2r^{2n}}\sqrt{\sum_{n=2}^\infty\frac{1}{n^2}}
\leq r+ r\left(\sqrt{\exp\left(\frac{4\lambda^2r^2}{1-r^2}\right)-1}\right)\!\!\left(\sqrt{\frac{\pi^2}{6}-1}\right).
$$
From \cite[Corollary 2.4]{Kim} we observe that $d(f(0), \partial f(\ID))\geq -F_\lambda(-1)$.
The inequality $(\ref{P3eq22})$
now holds whenever $r\leq r_0$ for some $r_0$, if we can show
$$
\phi(r):=r+r\sqrt{\exp\left(4\lambda^2r^2/(1-r^2)\right)-1}\sqrt{\pi^2/6-1}+ F_\lambda(-1)
$$
has one and only one root $r_0$ in $(0,1)$ and $\phi(r)\leq 0$ if and only if $r\leq r_0$.
It is easy to observe that $\phi(0)=F_\lambda(-1)<0$ and $\lim_{r\to1-}\phi(r)=\infty$ which together, by a use of intermediate value property
for continuous functions ensure the existence of one root $r_0\in(0,1)$ of $\phi(r)$. Again observing that
%$$
%\phi^\prime(r)=1+\sqrt{\frac{\pi^2}{6}-1}\left(\sqrt{\exp\left(\frac{4\lambda^2r^2}{1-r^2}\right)-1}
%+4\lambda^2r^2\frac{\exp\left(\frac{4\lambda^2r^2}{1-r^2}\right)}{(1-r^2)^2\sqrt{\exp\left(\frac{4\lambda^2r^2}{1-r^2}\right)-1}}\right)
%$$
$\phi^\prime(r)>0$ for all $r\in(0,1)$, we conclude that $\phi$ is strictly increasing in $(0,1)$.
This proves that $r_0$ is the only root of $\phi$ in $(0,1)$, and that
$\phi(r)\leq 0\Leftrightarrow r\leq r_0$.
\epf
\section{Bohr phenomenon for logarithmic power series}
In the first theorem of this section, we compute sharp Bohr radii for $\log(f(z)/z), z\in\ID$
and $\log(f^{-1}(w)/w), w\in\ID_{1/4}$ where $f\in\es(\mbox{or}\, \es^*)$. Besides, the sharp
Bohr radius for $\log(f(z)/z), f\in\mathcal{C}$ has been recorded in a subsequent remark.
\bthm\label{P3thm1}
Let $f\in\es(or\, \es^*)$ with $\log(f(z)/z)$ having Taylor expansion $(\ref{P3eq2})$. Then the inequality $(\ref{P3eq3})$
holds for $|z|=r\leq r_0=1-(1/\sqrt{e})=0.393\cdots$.
Moreover if the \textit{logarithmic coefficients} of $f^{-1}$ are given by $(\ref{P3eq11})$ where $f\in\es(or\, \es^*)$, then
inequality $(\ref{P3eq3})$ is satisfied for $|w|=r\leq r_0=(1/e)(\sqrt{e}-1)=0.238\cdots$.
Both results are sharp for the function $k_1(z)=z/(1+z)^2$.
%Also if $f\in\CC$, then the inequality $(\ref{P3eq3})$
%holds for $|z|=r\leq r_0=1-(1/e)=0.632\cdots$. This result is again sharp for the function $l(z)=z/(1-z)$.
\ethm
\bpf
For $f\in\es$, the following inequality is well known (see f.i. \cite[p. 722]{And}):
$$
\sum_{n=1}^\infty n|\gamma_n|^2r^n\leq \log\left(\frac{1}{1-r}\right)
$$
where $|z|=r$. Therefore a use of Cauchy-Schwarz inequality gives
$$
2\sum_{n=1}^\infty|\gamma_n|r^n\leq 2\sqrt{\sum_{n=1}^\infty n|\gamma_n|^2r^n}\sqrt{\sum_{n=1}^\infty\frac{r^n}{n}}
\leq 2\log\left(\frac{1}{1-r}\right),
$$
which is less or equal to 1 whenever $r\leq 1-(1/\sqrt{e})$.
Now to prove the second part of this theorem, we note that a use of the recent result \cite[Theorem 1]{Pon} gives
\be\label{P3eq12}
2\sum_{n=1}^\infty|\gamma_n|r^n\leq\sum_{n=1}^\infty\frac{1}{n}{2n\choose n}r^n.
\ee
It can be observed that for $r<1/4$,
$$
\sum_{n=1}^\infty{2n\choose n}r^n=\frac{1}{\sqrt{1-4r}}-1,
$$
or equivalently
$$
\sum_{n=1}^\infty{2n\choose n}r^{n-1}=\frac{4}{\sqrt{1-4r}(1+\sqrt{1-4r})}.
$$
Integrating both sides of the above equation from $0$ to $r$ we get
$$
\sum_{n=1}^\infty\frac{1}{n}{2n\choose n}r^n=\int_0^r\frac{4}{\sqrt{1-4x}(1+\sqrt{1-4x})}dx =:I.
$$
Setting $1-4x=t^2$, a little calculation reveals that
$$
I=2\int_{\sqrt{1-4r}}^1\frac{dt}{1+t}=2\log\frac{2}{1+\sqrt{1-4r}}.
$$
Therefore from $(\ref{P3eq12})$ it is seen that inequality $(\ref{P3eq3})$ will be satisfied whenever
$2\log(2/(1+\sqrt{1-4r}))\leq 1$, or, in other words, whenever $r\leq (1/e)(\sqrt{e}-1)$. Observing the fact that
the function $k_1(z)\in\es^*$, the sharpness of both the
results for the classes $\es$ and $\es^*$ can be shown from direct computations.
\epf

\br
For $f\in\CC$ with $\log(f(z)/z)$ having Taylor expansion $(\ref{P3eq2})$,
it is easy to prove the bounds $|\gamma_n|\leq 1/2n$ for $n\geq1$.
As a result the inequality $(\ref{P3eq3})$
holds for $|z|=r\leq r_0=1-(1/e)=0.632\cdots$. This result is sharp for the function $l(z)=z/(1-z)$.
\er

The next result includes Bohr phenomenon for $\log(f(z)/z)$ where $f\in\U(\lambda)$.
\bthm\label{P3thm3}
Suppose $0<\lambda\leq 1$ and $f\in\U(\lambda)$ with $\log(f(z)/z)$ having Taylor expansion $(\ref{P3eq2})$.
Then the inequality $(\ref{P3eq3})$ holds for
$$
|z|=r\leq r_0=
\begin{cases}
((1+\lambda)-\sqrt{(1+\lambda)^2-4\lambda(1-(1/e))})/2\lambda\,\,\mbox{if}\,\, \lambda\geq \lambda_0,\\
(1+\lambda^2)/2(1+\lambda)\,\,\mbox{if}\,\, \lambda< \lambda_0,\\
\end{cases}
$$
where $\lambda_0\approx 0.750792$ is the only root in $(0,1)$ of the equation
\be\label{P3eq7}
\lambda^5-2\lambda^4-2\lambda^3-(4/e)\lambda^2+(5-(8/e))\lambda+(2-(4/e))=0.
\ee
When $\lambda\geq \lambda_0$, this result is sharp for the function $k_\lambda(z)=z/(1+z)(1+\lambda z)$.
\ethm
\bpf
In \cite[Theorem 4]{Obra4}, it was shown that for $f\in\U(\lambda)$
%\be\label{P3eq4}
$$
\frac{f(z)}{z}\prec \frac{1}{(1+z)(1+\lambda z)}\,, z\in\ID,
%\ee
$$
which yields
$$
\log(f(z)/z)\prec -\log(1-z)-\log(1-\lambda z),
$$
or equivalently, in terms of Taylor expansions:
\be\label{P3eq5}
2\sum_{n=1}^\infty\gamma_nz^n\prec\sum_{n=1}^\infty\left(\frac{1+\lambda^n}{n}\right)z^n.
\ee
An application of \cite[Theorem 6.3]{Dur} on $(\ref{P3eq5})$ gives
\be\label{P3eq6}
4\sum_{n=1}^\infty\left(\frac{n}{1+\lambda^n}\right)|\gamma_n|^2r^n\leq
\sum_{n=1}^\infty\left(\frac{n}{1+\lambda^n}\right)\left(\frac{1+\lambda^n}{n}\right)^2r^n
=\sum_{n=1}^\infty\left(\frac{1+\lambda^n}{n}\right)r^n,
\ee
whenever $(n/(1+\lambda^n))r^n\geq ((n+1)/(1+\lambda^{n+1}))r^{n+1}$ for all $n\geq 1$, i.e.
whenever $r\leq ((1+\lambda^{n+1})/(1+\lambda^n))(n/(n+1))$ for all $n\geq 1$.
We now observe that
$$
u_n:=(1+\lambda^{n+1})/(1+\lambda^n)\,\mbox{and}\,v_n:= n/(n+1)
$$
both are monotonically increasing sequences, which imply $u_n\geq u_1=(1+\lambda^2)/(1+\lambda)$
and $v_n\geq v_1=1/2$, $n\in\IN$. Therefore the inequality $(\ref{P3eq6})$ remains
valid for $r\leq(1+\lambda^2)/2(1+\lambda)$.
Now a use of Cauchy-Schwarz inequality gives
$$
2\sum_{n=1}^\infty|\gamma_n|r^n\leq \sqrt{4\sum_{n=1}^\infty\left(\frac{n}{1+\lambda^n}\right)|\gamma_n|^2r^n}
\sqrt{\sum_{n=1}^\infty\left(\frac{1+\lambda^n}{n}\right)r^n}.
$$
Using $(\ref{P3eq6})$ on the right hand side of the above inequality we obtain
\be\label{P3eq8}
2\sum_{n=1}^\infty|\gamma_n|r^n\leq\sum_{n=1}^\infty\left(\frac{1+\lambda^n}{n}\right)r^n=-\log(1-r)-\log(1-\lambda r)
\ee
for $r\leq(1+\lambda^2)/2(1+\lambda)$. It is readily seen that the inequality $(\ref{P3eq3})$ would be satisfied if
$-\log(1-r)-\log(1-\lambda r)\leq 1$, i.e. if $r\leq r_0=\min\{r_b, (1+\lambda^2)/2(1+\lambda)\}$
where $r_b$ be the smallest positive root of
\be\label{P3eq9}
h(\lambda):=r^2-(1+(1/\lambda))r+(1/\lambda)(1-(1/e))=0.
\ee
It is clear that $r_b=((1+\lambda)-\sqrt{(1+\lambda)^2-4\lambda(1-(1/e))})/2\lambda$. We observe that
$h(0)=(1/\lambda)(1-(1/e))>0,\, h(1)=-1/\lambda e<0$, which by intermediate value property imply that
$h$ has at least one root in $(0,1)$. Also since $h^\prime(r)=2r-(1+(1/\lambda))<0$, $h$ is strictly decreasing in $(0,1)$,
and therefore $h$ has exactly one root in $(0,1)$ which is clearly $r_b$. Now it needs to be proved that $r_0=r_b$
precisely when $\lambda\geq\lambda_0$.
In fact we have to show that $r_b\leq (1+\lambda^2)/2(1+\lambda)$ if and only if $\lambda\geq\lambda_0$. In other words we need to
establish that $h((1+\lambda^2)/2(1+\lambda))\leq0$ if and only if $\lambda\geq\lambda_0$, or equivalently,
after a little calculation:
\be\label{P3eq10}
\lambda\geq\lambda_0\Leftrightarrow g(\lambda):=\lambda^5-2\lambda^4-2\lambda^3-(4/e)\lambda^2+(5-(8/e))\lambda+(2-(4/e))\leq0.
\ee
Now $g(0)=2-(4/e)>0$, and $g(1)=4-(16/e)<0$ assert that $g$ has at least one root in $(0,1)$,
which we choose to be $\lambda_0$. Now it is sufficient to show that
$g(\lambda)>0$ for $\lambda<\lambda_0$ and $g(\lambda)\leq 0$ for $\lambda\geq\lambda_0$.
We see that
$$
g^\prime(\lambda)=5\lambda^4-8\lambda^3-6\lambda^2-(8/e)\lambda+(5-(8/e)),
$$
and therefore
$g^\prime(0)=5-(8/e)>0$ and $g^\prime(1)=-4-(16/e)<0$, which ensure that $g^\prime$ has at least one root $\mu_0\in(0,1)$.
Again we observe that
$g^{\prime\prime}(\lambda)=20\lambda^3-24\lambda^2-12\lambda-(8/e)$, and therefore $g^{\prime\prime}(0)=-8/e$.
Since for $\lambda\in(0,1)$
$$
g^{\prime\prime\prime}(\lambda)=60\lambda^2-48\lambda-12=12(\lambda-1)(5\lambda+1)<0,
$$
$g^{\prime\prime}$ is strictly decreasing in $(0,1)$
and hence $g^{\prime\prime}(\lambda)< g^{\prime\prime}(0)<0.$ This now asserts that $g^\prime$ is strictly decreasing in $(0,1)$. Therefore
$\mu_0$ is the only root of $g^\prime$ in $(0,1)$, $g^\prime(\lambda)>0$ for $\lambda<\mu_0$ and $g^\prime(\lambda)< 0$ for $\lambda> \mu_0$.
Now let if possible, $\mu_0>\lambda_0$. Then $g$ is strictly increasing in $(0, (\lambda_0+\mu_0)/2)$, and as a result $g(\lambda_0)>g(0)>0$ which is
contrary to our assumption. Therefore $\mu_0\leq\lambda_0$, which shows that $g$ is strictly increasing in $(0, \mu_0)$ and strictly decreasing in
$(\mu_0, \lambda_0)\cup[\lambda_0,1)$, $\mu_0$ being a local maximum of $g$.
Clearly for any $\lambda\in(0, \mu_0)$, $g(\lambda)>g(0)>0$; for $\lambda\in[\mu_0, \lambda_0)$,
$g(\lambda)> g(\lambda_0)=0$; and for $\lambda\in[\lambda_0,1]$, $g(\lambda)\leq g(\lambda_0)=0$. This validates our assertion
$(\ref{P3eq10})$. Using Mathematica it can be computed that $\lambda_0$ is approximately $0.750792$. The sharpness part for $\lambda\geq\lambda_0$
is immediate from the lines of our computation.
\epf

\end{document}